\documentstyle[amscd,amssymb, verbatim,12pt]{amsart}

\theoremstyle{plain}
\newtheorem{Thm}{Theorem}
\newtheorem{Prop}[Thm]{Proposition}
\newtheorem{Cor}[Thm]{Corollary}
\newtheorem{Lem}[Thm]{Lemma}

 \theoremstyle{definition}
\theoremstyle{remark}

\numberwithin{equation}{section}

\begin{document} 
 \title{Action integrals along closed isotopies in coadjoint orbits} 
    
 \author{ ANDR\'{E}S   VI\~{N}A}
\address{Departamento de F\'{i}sica. Universidad de Oviedo.   Avda Calvo
 Sotelo.     33007 Oviedo. Spain. } 
 \email{cuevas@@pinon.ccu.uniovi.es}
\thanks{The author was partially supported by DGES, grant DGES-96PB0538}

 \maketitle
\begin{abstract}  Let ${\cal O}$ be the orbit of $\eta\in{\frak g}^*$ under the 
coadjoint action of the compact Lie group $G$. 
We give two formulae for calculating the action integral along 
a closed Hamiltonian isotopy on ${\cal O}$. The first one expresses this action 
in terms of a particular character of the isotropy 
subgroup of $\eta$. In the second one is involved the 
character of an irreducible representation of $G$.
 
 \end{abstract}

   \bigskip
 
 \section {Introduction} \label{S:intro}

\medskip

 Let $(M,\omega)$ be a quantizable \cite{nW92} symplectic manifold.
We will denote by $\text{Ham}(M)$ the group of Hamiltonian symplectomorphisms
\cite{Mc-S} of   $(M,\omega)$. In this note we will consider
  loops $\{\psi_t\}_{t\in[0,1]}$ in $\text{Ham}(M)$ at $\text{id}$.
  Given $q\in M$, the loop $\psi$ generates the closed curve 
$\{\psi_t(q)\,|\, t\in[0,1]\}$ in $M$ which is homologous to zero \cite[page 334]{Mc-S}.
As $(M,\omega)$ is quantizable,  
  it makes sense to define the action integral ${\frak A}_{\psi}(q)$ along such a curve
as the element of ${\Bbb R}/{\Bbb Z}$ given by
the formula \cite{aW89} \cite{Mc-S}
\begin{equation}\label{actint}
{\frak A}_{\psi}(q)=\int_S \omega-\int_0^1 f_t(\psi_t(q)) dt + {\Bbb Z},
\end{equation}
where $S$ is any 
$2$-surface whose boundary is the curve $\{\psi_t(q)\}$, and
where $f_t$ a fixed time dependent Hamiltonian  associated to $\{\psi_t\}$.

Since $(M,\omega)$ is quantizable, 
  one can choose a prequantum bundle  $L$ on $M$, endowed
 with a connection $D$ \cite{nW92}. On the other hand, let $X_t$ be 
 the corresponding
Hamiltonian vector fields determined by $f_t$, then  one can construct the operator
${\cal P}_t:=-D_{X_t}-2\pi if_t$, which acts on the sections of $L$. The equation
$\Dot \tau_t={\cal P}_t(\tau_t)$ defines a ``transport" of the
 section $\tau_0\in C^{\infty}(L)$ along $\psi_t$.
This transport enjoys the following nice property: If $D_Y\tau_0=0$, with $Y$ a vector 
field on $M$, then $D_{Y_t}\tau_t=0$, for $Y_t=\psi_t(Y)$ (see \cite{aV01}). From this
fact one can prove  that $\tau_1$ 
and $\tau_0$ differ in a  constant factor $\kappa(\psi)$; that is,
$\tau_1=\kappa(\psi)\tau_0$. A direct calculation shows that
$\kappa(\psi)=\text{exp}(2\pi i{\frak A}_{\psi}(q))$, where $q$ is an arbitrary point of $M$
 \cite{aV01}. Consequently the expression (\ref{actint})
  is independent of $q$, and  it makes sense  to define the action integral along $\psi$
by (\ref{actint}).

The purpose of this note is to calculate the value of the invariant $\kappa(\psi)$
when the manifold $M$ is a coadjoint orbit \cite{aK76} of a compact Lie group.
 However in Section 2 we study a more
general situation. If a Lie group $G$ acts
on the   manifold $M$ by symplectomorphisms and there is a moment map
for this action, then each $A\in{\frak g}$ determines a vector field $X_A$ on $M$
and  the corresponding Hamiltonian  $f_A$. We can construct the
respective operator ${\cal P}_A$
on $C^{\infty}(L)$, so one has a representation ${\cal P}$ of the Lie algebra ${\frak g}$ 
 on $C^{\infty}(L)$. When this representation extends to an action $\rho$ of the 
group $G$, the prequantization is said to be $G$-invariant.
In this case we will prove that  the value of $\kappa(\psi)$ can be expressed in terms of $\rho$.
More precisely, if the isotopy $\psi_t$ is determined by vector fields of type $X_{A_t}$
we show that $\tau_t=\rho(h_t)\tau_0$, where $h_t$ is the solution to   Lax equation 
$\Dot h_t h_t^{-1}=A_t$.

Section 3 is concerned with the invariant $\kappa(\psi)$ 
for closed isotopies $\psi$ in a coadjoint orbit
of a compact Lie group $G$. We study the value of $\kappa(\psi)$, when
the isotopy is  defined by  vector
 fields of type $X_{A_t}$. Given  $\eta\in{\frak g}^*$, 
the orbit ${\cal O}_{\eta}$
  of $\eta$ admits  a $G$-invariant prequantization if the prequantum bundle is defined by a
 character 
$\Lambda$ of $G_{\eta}$, the
subgroup of isotropy of $\eta$. In this case we prove 
that $\kappa(\psi)=\Lambda(h_1)$, with $h_t$
the solution to the corresponding Lax equation (Theorem \ref{Tmcoador}).

 If $G$ is a semisimple group, the choice of a maximal torus $T$ contained in $G_{\eta}$
permits us to define a $G$-invariant complex structure on $G/G_{\eta}={\cal O}_{\eta}$.
This complex structure, in turn, determines a holomorphic structure on $L$.
When the prequantization is $G$-invariant, ${\cal P}$ defines also 
a representation $\rho$ of $G$ on the space $H^0(L)$ of holomorphic sections of $L$. 
When $G_{\eta}$ itself is a maximal torus,  
  the Borel-Weil theorem allows us to characterize $\rho$ in terms of its 
highest weight. We prove that  the invariant $\kappa(\psi)$ for the closed isotopy
considered above is equal to $\chi(\rho)(h_1)/\text{dim}\,\rho.$ This fact permits us 
to calculate $\kappa(\psi)$ using the Weyl's character formula. This stuff is 
considered in Section 4.

In Section 5 we check the results of  Sections 3 and 4 in two particular cases.
In the first one we calculate directly  the value of $\kappa(\psi)$ for a
closed isotopy $\psi$ in ${\Bbb CP}^1$; Theorem \ref{Tmcoador}
   and    Weyl's character formula applied to this example give the same
result as the direct calculation. 
In \cite{aV01} we determined the value $\kappa(\psi)$
for a closed Hamiltonian flow $\psi$ in $S^2$; here we recover this number by applying
Theorem \ref{Tmcoador} to this isotopy.

\bigskip

\section {$G$-invariant prequantum data.}

\medskip

Let $G$ be a compact Lie group which acts on the left on the symplectic 
manifold $(M,\omega)$ by symplectomorphisms. We assume that this action is
Hamiltonian, and that $\Phi:M\rightarrow {\frak g}^*$ is a map moment for this action.

Given $A\in{\frak g}^*$, we denote by $X_A$,  the vector field 
on $M$ generated by $A$. Then  $(d\Phi(Y))\cdot A=\omega(Y, X_A)$, 
for any vector field $Y$ on $M$. The ${\Bbb R}$-valued map $\Phi\cdot A$ will be denoted by $f_A$;
so
\begin{equation}\label{moment}
\iota_{X_A}\omega=-df_A \,\,\,\,\text{and}\,\,\,\, \{f_A,f_B\}=\omega(X_ B,X_ A)=f_{[A,B]}.
\end{equation}

As we said one assumes that $(M,\omega)$ is quantizable.   Let $L$ be
 a prequantum bundle, i.e. $L$ is a Hermitian line
bundle over $M$ with a connection $D$, 
 whose curvature is $-2\pi i\omega$, then  one can define the prequantization map
\cite{jS80}
\begin{equation}\label{aingm}
 A\in{\frak g}\mapsto {\cal P}_A=-D_{X_ A}-2\pi if_A\in \text{End}(C^{\infty}(L)).
\end{equation}

\begin{Prop} The map ${\cal P}$ is a Lie algebra homomorphism.
\end{Prop}
\begin{pf} Since the action of $G$ is on the left, the map
$A\in{\frak g}\mapsto X_ A\in\Xi(M)$, where $\Xi(M)$ denotes the set of 
vector fields on $M$, is a Lie Algebra antihomomorphism (see \cite{KN} p.42); that is, 
\begin{equation}\label{antiLi}
X_{[A,B]}=-[X_A,X_B]
\end{equation}

On the other hand, if $\tau$ is a section of $L$ 
\begin{equation}\label{calPA}
[{\cal P}_A, {\cal P}_B]\tau=[D_{X_A},D_{X_B}]\tau+4\pi i\omega(X_ A,X_ B)\tau.
\end{equation}
 Since the curvature of $D$ is $-2\pi i\omega$ 
$$-2\pi i\omega(X_A,X_B)\tau= [D_{X_A},D_{X_B}]\tau + D_{[X_A,X_B]}\tau.$$
Using (\ref{antiLi}),  
   (\ref{moment}) and (\ref{calPA}) one obtains
$$[{\cal P}_A,{\cal P}_B]\tau
={\cal P}_{[A,B]}\tau.$$
\end{pf}

The prequantum data $(L,D)$ are said to be $G$-invariant, if there is a action 
$\rho$ of $G$ on $C^{\infty}(L)$ which generates ${\cal P}$ \cite{vG82}. Henceforth in this Section
we assume that 
the prequantum data are $G$-invariant.

Let $\{A_t \}_t$ be a curve in ${\frak g}$ with $A_0=0$. Given $\tau\in C^{\infty}(L)$
we consider the equation for the section $\tau_t$ of $L$
\begin{equation}\label{transport}
\frac{d\tau_t}{dt}={\cal P}_{A_t}(\tau_t),\,\,\,\,\, \tau_0=\tau
\end{equation}
This is the equation of the ``transport" of the section $\tau$ along the isotopy determined 
by the vector fields $X_{A_t}$ (see \cite{aV01}). We will try to find a curve $h_t$ in $G$,
such that $h_0=e$ and $\rho(h_t)(\tau)=\tau_t,$ where $\tau_t$ is solution to
(\ref{transport}). As $\rho:G\rightarrow \text{Diff}(C^{\infty}(L))$ is a group homomorphism,
 $\rho\circ {\cal L}_g={\cal L}_{\rho(g)}\circ\rho$,
where ${\cal L}_a$ is the left multiplication by $a$ in the respective group.
 The corresponding tangent maps satisfy 
\begin{equation}\label{estrella}
\rho_*\circ {{\cal L}_g}_*={{\cal L}_{\rho(g)}}_*\circ\rho_*.
\end{equation}
 If we put $F_t$ for  diffeomorphism   $\rho(h_t)=:F_t$, and we define $Y_t\in{\frak g}$ by
$$\Dot h(t)={{\cal L}_{h(t)}}_*(Y_t),$$
then by (\ref{estrella})
\begin{equation}\label{aux1}
\frac{d F_t}{dt}=\rho_*(\Dot h_t)={{\cal L}_{F_t}}_*({\cal P}(Y_t)).
\end{equation}
As ${\cal L}_{F_t}(C)=F_t\circ C$, if $C\in\text{End}(C^{\infty}(L))\subset\Xi(C^{\infty}(L))$, then
(\ref{aux1}) can be written
$$\frac{d F_t}{dt}=F_t\circ {\cal P}(Y_t).$$
If we introduce this formula  in (\ref{transport}), we obtain
$$\frac{d\tau_t}{dt}=(F_t\circ {\cal P}(Y_t))\tau=({\cal P}_{A_t}\circ F_t)\tau.$$
Hence
\begin{equation}\label{F_t}
F_t\circ {\cal P}_{Y_t}\circ F_t ^{-1}={\cal P}_{A_t}
\end{equation}
Let $\{m(u)\}_u$ a  curve in  $G$ which defines $Y_t\in{\frak g}$, then
$$F_t\circ{\cal P}_{Y_t}\circ F_t^{-1}=
\frac{d}{du}\biggr|_{u=0}\rho(h_tm(u)h_t^{-1}).$$
By (\ref{F_t}) one can take $Y_t=\text{Ad}_{h_t^{-1}}A_t$; so $h_t$ is the solution to the 
Lax equation
\begin{equation}\label{Lax}
\Dot h_th_t^{-1}=A_t\,\,\,\,\,\,\,\, h_0=e.
\end{equation}
We have proved
\begin{Thm}\label{ThmA} The solution $\tau_t$ to (\ref{transport}) is given by $\rho(h_t)\tau$,
where $h_t$ satisfies  equation (\ref{Lax}).
\end{Thm}

\smallskip

Let $\{A_t\,|\, t\in[0,1]  \}$ be a curve in ${\frak g}$ such that the Hamiltonian isotopy
$\{\psi_t \}_{t\in[0,1]}$ generated by the vector fields $X_{A_t}$ is closed; i.e.
$\psi_0=\psi_1=\text{id}$. We have proved in \cite{aV01} that 
 if $\tau_t$ is the solution of (\ref{transport}), then
\begin{equation}\label{kappa}
\tau_1=\kappa(\psi)\tau,
\end{equation}
 for every $\tau\in C^{\infty}(L)$, where 
$\kappa(\psi)=\text{exp}\big(2\pi i{\frak A}_{\psi}(q)  \big)$,
and ${\frak A}_{\psi}(q)$ is the action integral along the curve 
$\{\psi_t(q) \}_t$, for $q$ arbitrary in $M$.  
On the other hand, if $h_t$ is a curve in $G$ solution to (\ref{Lax}), by Theorem  
\ref{ThmA} $\tau_1=\rho(h_1)(\tau)$. It follows from (\ref{kappa}) that
$\rho(h_1)=\kappa(\psi)\,\text{Id}.$ Thus we have
  
\begin{Cor}\label{Corchi}If $W$ is a  finite dimensional $\rho$-invariant subspace of $C^{\infty}(L)$,
and $\rho_W$ is the restriction of $\rho$ to this subspace, then for the character of 
$\rho_W$ holds the following formula
$$\chi(\rho_W)(h_1)=\kappa(\psi)\,\text{dim}(W).$$
\end{Cor}

\bigskip

\section {The invariant $\kappa(\psi)$ in a coadjoint orbit} 

\medskip
 
Let $G$ be a compact Lie group, and we consider the coadjoint action of $G$
 on ${\frak g}^*$ defined by 
$$(g\cdot\eta)(A)=\eta(g^{-1}\cdot A),$$
for $g\in G$, $\eta\in{\frak g}^*$, $A\in{\frak g}$ and $g\cdot A=\text{Ad}_g A$
(see \cite{aK76} \cite{nW92}).

If $X_A$ is the vector field on ${\frak g}^*$ determined by $A$, the map
 $l_g:\mu\in{\frak g}^*\mapsto g\cdot \mu\in{\frak g}^*$
satisfies
\begin{equation}\label{l_g_*}
(l_g)_*(X_A(\mu))=X_{g\cdot A}(g\cdot \mu).
\end{equation}

Given $\eta\in{\frak g}^*$, by ${\cal O}_{\eta}=:{\cal O}$ will be denoted the orbit of $\eta$
under the coadjoint action of $G$.   
 On ${\cal O}$
one can consider the $2$-form $\omega$ determined by
\begin{equation}\label{omeganu}
 \omega_{\nu}(X_{A}(\nu), X_{B}(\nu))=\nu([A,B]).
\end{equation}
This $2$-form defines a symplectic structure on ${\cal O}$, and the action of $G$
preserves $\omega$. For each $A\in{\frak g}$ one defines 
the function $h_A\in C^{\infty}({\cal O})$ by $h_A(\nu)=\nu(A)$, and for this function holds the 
formula
\begin{equation}\label{iotah}
\iota_{X_{A}}\omega=dh_A.
\end{equation}

The orbit ${\cal O}$ can be identified with $G/G_{\eta}$, where $G_{\eta}$
is the subgroup of isotropy of
 $\eta$. The Lie algebra of this subgroup is    
$${\frak g}_{\eta}=
\{A\in{\frak g}\,|\, \eta([A,B])=0, \,\,\text{for every}\,\,B\in{\frak g}  \}$$

The orbit  ${\cal O}$ possesses a $G$-invariant prequantization iff the linear functional
\begin{equation}\label{intfunc}
\lambda:C\in{\frak g}_{\eta}\mapsto 2\pi i\eta(C)\in i{\Bbb R}
\end{equation}
is integral; i. e., iff there is a character $\Lambda:G_{\eta}\rightarrow U(1)$
whose derivative is the functional (\ref{intfunc}) (see \cite{bK70}).
Henceforth we assume the existence of such a character $\Lambda$. The corresponding
prequantum bundle $L$ over ${\cal O}=G/G_{\eta}$ is defined by
$L=G\times_{\Lambda}{\Bbb C}=(G\times {\Bbb C})/\simeq$, with
$(g,z)\simeq (gb^{-1}, \Lambda(b)z)$, for $b\in G_{\eta}$.

Each section $\sigma$ of $L$ determines a $\Lambda$-equivariant function $s:G\rightarrow{\Bbb C}$
by the relation
\begin{equation}\label{sfunction}
\sigma(gG_{\eta})=[g,s(g)].
\end{equation}

 The ${\Bbb C}^{\times}$-principal bundle associated to 
$L$ is $L^{\times}=L-\{\text{zero section}\}$. 
If $v$ denotes the element $[e,1]\in L^{\times}$, then 
$T_v(L^{\times})\simeq({\frak g}\oplus {\Bbb C})/{\frak f}_v$,
with
$${\frak f}_v=\{(B,\,-2\pi i\eta(B)\,|\, B\in{\frak g}_{\eta}  \}.$$
The connection form $\Omega$ on $L^{\times}$ is constructed in \cite{bK70} p.198.  
The form $\Omega$ can be written $\Omega=(\theta, \,d)$, where $\theta$ is the left 
invariant form on $G$ whose value at $e$ is $\eta$, and 
$d\in\text{Hom}_{\Bbb C}({\Bbb C},{\Bbb C})$ is defined by $d(z)=(2\pi i)^{-1}z$.
 It is clear that 
$\Omega_v$ vanishes on ${\frak f}_v$ and that it defines an element of $T_v^*(L^{\times})$.

On the other hand the section 
 $\sigma$ determines a lift
$\sigma^{\sharp}:L^{\times}\rightarrow {\Bbb C}$ by the formula
\begin{equation}\label{sostenido}
\sigma(\pi(y))=\sigma^{\sharp}(y)y,
\end{equation}
here $\pi:L\rightarrow {\cal O}$ is the projection map. It follows from (\ref{sfunction}) 
\begin{equation}\label{ssostenido}
s(g)=\sigma^{\sharp}([g,z])z.
\end{equation}

 We denote by ${\cal E}_{\Lambda}$ the space of $\Lambda$-equivariant functions on $G$.
The identification $C^{\infty}(L)\thicksim {\cal E}_{\Lambda}$ allows us to translate the 
action ${\cal P}$ defined in (\ref{aingm}) to a
 representation of ${\frak g}$ on ${\cal E}_{\Lambda}$.

\begin{Thm}\label{thms} The action ${\cal P}$ on ${\cal E}_{\Lambda}$ 
 is given by ${\cal P}_A(s)=-R_A(s)$,
where $R_A$ is the right invariant vector field on $G$ determined by $A$.
\end{Thm}
\begin{pf}
Let $\sigma$ be a section of $L$, by (\ref{iotah}) 
${\cal P}_A(\sigma)=-D_{X_A}\sigma+2\pi ih_A\sigma.$ 
We will determine the lift $({\cal P}_A(\sigma))^{\sharp}$.

 The vector $X_{A}(g\cdot\eta)\in T_{g\cdot\eta}({\cal O})$ is defined by the curve
$u\mapsto e^{uA}g\cdot\eta$ in ${\cal O}$. A lift of this curve at the point $[g,z]\in L^{\times}$
will be a curve of the form
$\gamma(u)=[e^{uA}g,\, z_u]$, with $z_u=ze^{ux}$. The vector tangent to $\gamma$ at $[g,z]$ 
is $\Dot\gamma(0)=[R_{A}(g),x]$, where $R_{A}(g)$ is 
the value at $g$ of the right invariant vector field in $G$ defined by $A$.

The condition $\Omega(\Dot\gamma(0))=0$ implies 
\begin{equation}\label{x0-2}
x=-2\pi i\eta(g^{-1}\cdot A).
\end{equation}
Therefore the horizontal lift of $X_{A}(g\cdot\eta)$ is
$$X_{A}^{\sharp}([g,z])=[R_{A}(g),\, -2\pi i\eta(g^{-1}\cdot A)],$$
and by  (\ref{ssostenido}) the action of $X_{A}^{\sharp}([g,z])$ on 
the function $\sigma^{\sharp}$
can expressed in terms of $s$
$$X_{A}^{\sharp}([g,z])(\sigma^{\sharp})=
\frac{d}{du}\biggr|_{u=0}\Big( \frac{s(e^{uA}g)}{ze^{ux}} \Big)=\frac{R_{A}(g)(s)}{z}-
\frac{xs(g)}{z}.$$
Since $X_{A}^{\sharp}(\sigma^{\sharp})=(D_{X_A}\sigma)^{\sharp}$ \cite[page 115]{KN}, from
(\ref{x0-2})
and (\ref{ssostenido}) it turns out that 
 the equivariant function associated to $D_{X_{A}}\sigma$ is 
\begin{equation}\label{equivarDat}
g\in G\mapsto R_{A}(g)(s)+2\pi i\eta(g^{-1}\cdot A)s(g)\in{\Bbb C}.
\end{equation}
 
Obviously the equivariant function defined by the section $h_{A}\sigma$ is
the function $\lambda_{A}s$, where
 $\lambda_{A}(g)=h_{A}(gG_{\eta})=(g\cdot\eta)(A)=\eta(g^{-1}\cdot A).$
It follows from (\ref{equivarDat})   that
the equivariant function which corresponds to
$-D_{X_A}\sigma+2\pi ih_A\sigma$ is $-R_A(s)$.
\end{pf}

\begin{Cor}\label{CorII} The action ${\cal P}$ on ${\cal E}_{\Lambda}$ is induced by the action
$$\rho:(b,s)\in G\times {\cal E}_{\Lambda}\mapsto s\circ{\cal L}_{b^{-1}}\in{\cal E}_{\Lambda},$$
where ${\cal L}_c$ is left multiplication by $c$   in the group $G$.
 \end{Cor}
\begin{pf}If $g_t=e^{tA}\in G$, then
$$\frac{d\,\rho_{g_t}(s)}{dt}\biggr|_{t=0}(g)=\frac{d}{dt}\biggr|_{t=0}s(e^{-tA}g)=-R_A(g)(s)={\cal P}_A(s)(g).$$
\end{pf} 

From Corollary \ref{CorII} it follows that the prequantum data $(L,D)$  are $G$-invariant.

\smallskip

Let $\{\psi_t\,|\, t\in[0,1] \}$ be a closed Hamiltonian isotopy on ${\cal O}$; that is,
 a Hamiltonian isotopy such that $\psi_1=\text{id}$. We also assume that
 the corresponding Hamiltonian vector fields
are invariant; that is,
$$\frac{d\psi_t(q)}{dt}=X_{A_t}(\psi_t(q)),\,\,\,\text{with}\,\,\,A_t\in{\frak g}.$$

If $\sigma$ is a section of $L$,  $\sigma_t$ will denote the solution to the equation
\begin{equation}\label{transportt}
\frac{d\sigma_t}{dt}= {\cal P}_{A_t}(\sigma_t),\,\,\,\,\,\,\, \sigma_0=\sigma.
\end{equation}
 By Theorem \ref{thms}, equation (\ref{transportt}) on the points
$\{g_t\}_{t\in[0,1]}$ of a curve in $G$ gives rise to   
\begin{equation}\label{Dotst}
\Dot s_t(g_t)=-R_{A_t}(g_t)(s_t),
\end{equation}
for the corresponding equivariant functions.
In particular, if $g_t$ is the curve such that $g_0=e$ and 
$\Dot g_t=R_{A_t}(g_t)\in T_{g_t}(G);$ in other words, $g_t$ satisfies the Lax equation
    $\Dot g_tg_t^{-1}=A_t$, then
$$R_{A_t}(g_t)(s_t)=\frac{d}{du}\biggr|_{u=t}s_t(g_u).$$
Using (\ref{Dotst}) one deduces
\begin{equation}\label{constanc}
\Dot s_t(g_t)+\Dot g_t(s_t)=0
\end{equation}
If we consider the function $w:[0,1]\rightarrow{\Bbb C}$ defined by
$w_t=s_t(g_t)$; by
(\ref{constanc}) $w$ is constant. So $s_1(g_1)=s_0(e).$
If $g_1\in G_{\eta}$, as $s_1$ is $\Lambda$-equivariant
$s_1(g_1)=\Lambda(g_1^{-1})s_1(e)$; so
\begin{equation}\label{sigma1}
\sigma_1(eG_{\eta})=\Lambda(g_1)\sigma_0(eG_{\eta}).
\end{equation} 
The following Theorem is consequence of (\ref{kappa}) and (\ref{sigma1})
\begin{Thm}\label{Tmcoador} If $\{\psi_t \}$ is the closed 
Hamiltonian isotopy in ${\cal O}$ generated by the 
vector fields $\{X_ {A_t} \}$, then $\kappa(\psi)=\Lambda(g_1)$, where $g_t\in G$
is the solution to $\Dot g_tg_t^{-1}=A_t$, with $g_0=e$ and $g_1\in G_{\eta}$.
\end{Thm} 

{\sc Remark.} Theorem \ref{Tmcoador} can also be deduced as a consequence of Theorem \ref{ThmA} 
and Corollary \ref{CorII}. In fact, if $h_t$ is the solution to
 $\Dot h_th_t^{-1}=A_t$, with the introduced notations
$$\sigma_1(a)=(\rho(h_1)\sigma)(a)=[a, s(h_1^{-1}a)]=\Lambda(h_1)[a,s(a)]=\Lambda(h_1)\sigma(a).$$

\bigskip

\section {Relation with Weyl's character formula.}

Let us assume that $G$ is semisimple Lie group \cite{wF91}, and let $T$ a maximal torus
with $T\subset G_{\eta}$ (see \cite{vG82} p.166). One has the corresponding 
decomposition of ${\frak g}_{\Bbb C}={\frak g}\otimes_{\Bbb R}{\Bbb C}$
in direct sum of root spaces
$${\frak g}_{\Bbb C}={\frak h}\oplus\sum {\frak g}_{\alpha},$$
where ${\frak h}={\frak t}_{\Bbb C}$, and $\alpha$ ranges over the set of roots.
This decomposition gives the real counterpart
$${\frak g}={\frak t}\oplus\sum_{\alpha\in P}
\big({\frak g}_{\alpha}\oplus {\frak g}_{-\alpha}\big)\cap{\frak g},$$
where $P$ is a set of positive roots. We denote by $\alpha^{\vee}$ the element
of $[ {\frak g}_{\alpha} , {\frak g}_{-\alpha}]$ such that $\alpha(\alpha^{\vee})=2$,
while $\beta(\alpha^{\vee})\in{\Bbb Z}$ for every root $\beta$. 

$\eta$ extends in a natural way to ${\frak g}_{\Bbb C}$. If $Y\in{\frak g}_{\alpha}$, then 
as $\alpha^{\vee}\in{\frak g}_{\eta}$
$$0=\eta([\alpha^{\vee},Y])=2\eta(Y).$$
Hence $\eta$ vanishes on $\sum {\frak g}_{\alpha}.$ If $\eta(\alpha^{\vee})\ne 0$, for all
root $\alpha$, then ${\frak g}_{\eta}={\frak t}$; in this case $\eta$ is said to be regular.
Henceforth we assume that $\eta$ is regular. 

We define ${\frak b}={\frak h}\oplus {\frak n}$, where 
$${\frak n}=\sum_{\alpha\in P}{\frak g}_{\alpha}.$$
Then ${\frak b}$ is a Borel subalgebra of ${\frak g}_{\Bbb C}$, which corresponds
to a Borel subgroup $B$ of $G$.

We have 
$$T_{\eta}({\cal O})={\frak g}/{\frak g}_{\eta}=
\sum_{\alpha\in P}\Big( {\frak g}_{\alpha}\oplus{\frak g}_{-\alpha} \Big)\cap{\frak g}.$$
Hence 
$$T^{\Bbb C}_{\eta}({\cal O})=\sum_{\alpha\in P}\Big(  {\frak g}_{\alpha}\oplus{\frak g}_{-\alpha} \Big ).$$
One defines 
$$T^{0,1}_{\eta}{\cal O}:={\frak n},$$
and 
$$T^{0,1}_{g\cdot\eta}{\cal O}:=\{ X_{g\cdot A}(g\cdot\eta)\,|\,A\in{\frak n}\}.$$
If $g_1\cdot\eta=g_2\cdot \eta$, then $g_1^{-1}g_2\in T$.
As ${\frak g}_{\alpha}$ is an eigenspace for the action of $T$, then 
$g_1^{-1}g_2\cdot A\in{\frak n}$, if $A\in{\frak n}$. Therefore the
spaces $T^{0,1}_{g\cdot\eta}$ are well-defined.

 For  
$A\in{\frak n}$, one can define the vector field ${\cal A}$ on  ${\cal O}$
by ${\cal A}(g\cdot\eta)=X_{g\cdot A}(g\cdot \eta)$.
By (\ref{l_g_*}) $(l_g)_*{\cal A}={\cal A}$, hence the above complex foliation
defined on ${\cal O}$ is $G$-invariant. Since
the vector $X_{g\cdot A}(g\cdot\eta)$ is defined by the curve
$e^{tg\cdot A}g\cdot \eta=ge^{tA}\cdot \eta$, then     the left invariant vector field $L_A$
on $G/T$ is the field which
 corresponds to  ${\cal A}$,
  in the identification of $G/T$ with ${\cal O}$.

The vector spaces $T^{1,0}$ are defined in the obvious way. As ${\frak n}$
is a subalgebra of ${\frak g}_{\Bbb C}$, the decomposition
$T^{\Bbb C}({\cal O})= T^{1,0}\oplus T^{0,1}$ define a complex structure on 
${\cal O}$. This complex manifold can be identified with $G_{\Bbb C}/B$.

We assume that the integral
 functional $\lambda$ in (\ref{intfunc}) satisfies   $\lambda(\alpha^{\vee})\leq 0$
for every $\alpha\in P$; this means that $-\lambda$ is a dominant
 weight for $T$ \cite{JJ00}.
Using the complex structure on ${\cal O}=G/T$ and the covariant derivative $D$
 on the prequantum bundle $L=G\times_{\Lambda}{\Bbb C}$, it is possible to
define a holomorphic structure in $L$. The section $\tau$ of $L$
is said to be holomorphic iff $D_Z\tau=0$ for any vector field $Z$
of type $(0,1)$. In this way $L$ can be regarded as a holomorphic line bundle
over $G_{\Bbb C}/B$. The homomorphism $\Lambda:T\rightarrow U(1)$
extends trivially to $B$, since $B$ is a semidirect product of $H=T_{\Bbb C}$ and
the nilpotent subgroup whose Lie algebra is ${\frak n}$; and  each   section $\sigma$
of $L$ determines a  function $s:G_{\Bbb C}\rightarrow {\Bbb C}$ which
is $\Lambda$-equivariant. 

Given $A\in {\frak n}$, the Proof of Theorem \ref{thms} shows that 
the equivariant function associated
to $D_{\cal A}\sigma$ is the map
$$g\in G_{\Bbb C}\mapsto R_{g\cdot A}(g)(s)+2\pi i\eta(g^{-1}g\cdot A)s(g)\in{\Bbb C}.$$
As $\eta$ vanishes on $\frak n$ and $R_{g\cdot A}(g)=L_A(g)$, the function associated
to $D_{\cal A}\sigma$ is $L_A(s)$. The section $\sigma$ is holomorphic
if $D_{\cal A}\sigma=0$, for every $A\in{\frak n}$; in this case  
$L_A(s)=0$ for $A\in{\frak n}$, that is, 
 $s$ is a holomorphic function on $G_{\Bbb C}$. So 
  the space $H^0(G_{\Bbb C}/B, L)$ is isomorphic to the space
$${\cal E}_{\Lambda,P}:=\{s:G_{\Bbb C}\rightarrow {\Bbb C}\,|
\, s \,\,\text{is holomorphic and}\,\,\,\Lambda-\text{equivariant}   \}.$$

The Borel-Weil Theorem asserts that the action of $G$ on 
the space ${\cal E}_{\Lambda,P}$ given by 
$g\star s=s\circ{\cal L}_{g^{-1}}$ is an irreducible representation of $G$;
more precisely the contragredient representation of that one whose highest weight
is $-\lambda$ (see \cite{JJ00} pages 290, 300).

\begin{Lem}\label{LemAin} If $A\in{\frak n}$,  then $[{\cal A},X_B]=0$ for any
$B\in{\frak g}_{\Bbb C}$. 
\end{Lem}
 \begin{pf}   The flow $\varphi_t$
determined by $X_B$ is given $\varphi_t(g\cdot\eta)=e^{tB}g\cdot\eta$.   And the flow $\phi_t$
of ${\cal A}$ is $\phi_t(g\cdot\eta)=e^{tg\cdot A}g\cdot\eta=ge^{tA}\cdot\eta$.
Hence 
$$(\varphi_t\circ\phi_t)(g\cdot\eta)=e^{tB}ge^{tA}\cdot\eta=
(\phi_t\circ\varphi_t)(g\cdot\eta).$$
\end{pf}
 
\begin{Prop}\label{PropDcalA} If $D_{\cal A}\sigma=0$ for any $A\in{\frak n}$, then
 $D_{\cal A}{\cal P}_B\sigma=0$
for any $B\in{\frak g}$.
\end{Prop}

\begin{pf}  Since  $D_{\cal A}\sigma=0$, it follows from (\ref{aingm})
\begin{equation}\label{DcalA}
D_{\cal A}({\cal P}_B\sigma)=-D_{\cal A}D_{X_B}\sigma
+2\pi i{\cal A}(h_B)\sigma.
\end{equation}
As 
$$[D_{\cal A},D_{X_B}]\sigma=D_{[{\cal A},X_B]}\sigma-2\pi i\omega({\cal A},X_B)\sigma,$$
from (\ref{DcalA}) and (\ref{iotah}) we deduce
$$D_{\cal A}({\cal P}_B\sigma)=-D_{[{\cal A},X_B]}\sigma.$$
Now the proposition is consequence of Lemma \ref{LemAin}. 
\end{pf}
 A direct consequence of Proposition \ref{PropDcalA} is
\begin{Cor}${\cal P}$ defines a representation of ${\frak g}$ on $H^0(G_{\Bbb C}/B, L)$.
 \end{Cor}
Denoting by $\pi$ the irreducible representation of $G$
whose highest weight is $-2\pi i\eta$, and by $\pi^*$ its
dual, we have
\begin{Cor} The representation ${\cal P}$ on $H^0(G_{\Bbb  C}/B, L)$ is the 
derivative    of $\pi^*$.
\end{Cor}
 \begin{pf} It is a consequence of Corollary \ref{CorII} and Borel-Weil theorem
\end{pf}

The subspace ${\cal E}_{\Lambda,P}\subset{\cal E}_{\Lambda}$ is invariant under the representation
$\rho$ defined in Corollary \ref{CorII}, and the restriction of $\rho$ to
${\cal E}_{\Lambda,P}$ is precisely the representation $\pi^*$.
  From by Corollary \ref{Corchi}
it follows
\begin{Thm} Let $\eta$ be an element of ${\frak g}^*$, such that
the orbit ${\cal O}_{\eta}$  is quantizable and $-2\pi i\eta$ is a dominant weight
 for the maximal torus
$G_{\eta}$. 
If $\{\psi_t \}$ is the closed Hamiltonian 
isotopy in ${\cal O}_{\eta}$ generated by the 
vector fields $\{X_ {A_t} \}$, then
\begin{equation}\label{Weylf}
\kappa(\psi)=\frac{\chi(\pi^*)(h_1)}{\text{dim}\,\pi},
\end{equation}
 where $h_t\in G$
is the solution to $\Dot h_th_t^{-1}=A_t,\,h_0=e$, and  $\pi$ is the representation of
$G$ whose highest weight is $-2\pi i\eta.$
\end{Thm} 
Now the character $\chi(\pi^*)$ and the dimension $\text{dim}\,\pi$ can be determined by
  Weyl's character formula \cite{JJ00}, and so $\kappa(\psi)$.

\bigskip

\section {Examples}

{\sc The invariant $\kappa(\psi)$ in ${\Bbb CP}^1$.}
Let $G$ be the group $SU(2)$ and $\eta$ the element of ${\frak su}(2)^*$
defined by
\begin{equation}\label{eta}
\eta\begin{pmatrix}
ci& w \\
-\bar w& -ci
\end{pmatrix}=-\frac{c}{2\pi}.
\end{equation}
The subgroup of isotropy $G_{\eta}$ is $U(1)\subset SU(2)$, so
 the coadjoint orbit ${\cal O}_{\eta}$ can be identified with 
$SU(2)/U(1)={\Bbb CP}^1$.
The element
$$g=\begin{pmatrix} z_0 & -\bar z_1 \\
                     z_1 & \bar z_0
\end{pmatrix}\in SU(2)$$ determines the point
$(z_0:z_1)\in{\Bbb CP}^1$. Hence to $\eta\in{\cal O}_{\eta}$ corresponds
$p=(1:0)\in{\Bbb CP}^1$. 
 For $z_0\ne 0$ we put $(z_0:z_1)=(1:x+iy)$.

Denoting by $A$ and $B$ the following matrices of ${\frak su}(2)$
\begin{equation}\label{matrices}
 A:=\begin{pmatrix} 0 & i \\
                     i & 0 
\end{pmatrix}, \,\,\,
B:=\begin{pmatrix} 0 & 1 \\
                     -1 & 0 
\end{pmatrix},
\end{equation}
by (\ref{omeganu}) 
\begin{equation}\label{omegaeta}
\omega_{\eta}(X_A,X_B)=\eta([A,B])=\frac{1}{\pi}.
\end{equation} 
As 
$$e^{tA}=\begin{pmatrix} \cos t & i\sin t \\
                         i\sin t & \cos t
\end{pmatrix}, $$
  the curve  $\{e^{tA}\eta\}$,  which defines $X_ A(p)$, is $(\cos t :i\sin t)$. 
Hence $X_A (p)$, 
  expressed in the real coordinates $(x,y)$, is equal to
 $\big(\frac{\partial}{\partial y}\big)_p$. Similarly
$X_B(p)= -\big(\frac{\partial}{\partial x}\big)_p$. Hence it follows from (\ref{omegaeta})
$$\omega_p=\frac{1}{\pi}dx\wedge dy=\frac{i}{2\pi}dz\wedge d\bar z,$$
where $z=x+iy$. Therefore $({\cal O}_{\eta},\omega)$ can be identified with ${\Bbb CP}^1$
endowed with the Fubini-Study form
\begin{equation}\label{Fubini-Study}
\omega=\frac{i}{2\pi}\frac{dz\wedge d\bar z}{(1+z\bar z)^2}=
\frac{1}{\pi}\frac{1}{(x^2+y^2+1)^2}dx\wedge dy.
\end{equation}

Let us consider the symplectomorphism $\psi_t$ on ${\Bbb CP}^1$ defined by
$$(z_0:z_1)\in{\Bbb CP}^1\mapsto (e^{-ia_t}z_0:e^{ia_t}z_1)\in{\Bbb CP}^1,$$
where $a_t\in{\Bbb R}$. If we assume that $a_0=0$ and $a_1=k\pi$, with $k\in{\Bbb Z}$,
then $\{ \psi_t \,|\, t\in[0,1]\}$ is a closed Hamiltonian isotopy on ${\Bbb CP}^1$.
We will determine $\kappa(\psi)$ by direct calculation. In real coordinates
\begin{equation}\label{psit}
\psi_t(x,y)=\big(x\cos 2a_t-y\sin 2a_t,\,x\sin 2a_t +y\cos 2a_t            \big).
\end{equation}
 A straightforward calculation shows that the Hamiltonian vector field $X_t$ defined
by
$$\frac{d\psi_t(q)}{dt}=X_t(\psi_t(q))$$
is $X_t(x,y)=2\Dot a_t\big(-y\frac{\partial}{\partial x}+x\frac{\partial}{\partial y} \big).$
It follows from (\ref{Fubini-Study})
$$\iota_{X_t}\omega=\frac{-2\Dot a_t}{\pi(x^2+y^2+1)^2}\big(xdx+ydy   \big).$$
A Hamiltonian function $f_t$ associated to $X_t$ is 
$$f_t(x,y)=-\frac{\Dot a_t}{\pi(x^2+y^2+1)}+c_t,$$
$c_t$ being a constant. If we impose $\int f_t\omega=0$, then
$$c_t=c_t \int_{{\Bbb CP}^ 1} \omega=
\frac{\Dot a_t}{\pi^2}\int_{{\Bbb CP}^ 1}\frac{1}{(x^2+y^2+1 )^3}dx\wedge dy=\frac{\Dot a_t}{2\pi}.$$
Thus the normalized Hamiltonian function is 
$$f_t(x,y)=\frac{\Dot a_t}{2\pi}\Big(\frac{x^2+y^2-1}{x^2+y^2+1}\Big).$$

Given $q=(x_0,y_0)\in{\Bbb CP}^1$, from (\ref{psit}) it follows that the
set 
$$\{\psi_t(x_0,y_0)\,|\,t\in[0,1]\}$$
 is a circle  in the plane $(x,y)$ with centre at 
$(0,0)$; therefore
\begin{equation}\label{auxr1}
\int_0^1f_t(\psi_t(q))dt=\frac{k}{2}\Big(\frac{x_0^2+y_0^2-1}{x_0^2+y_0^2+1}  \Big).
\end{equation}
On the other hand the $1$-form $\theta=(-x+iy)(x^2+y^2+1)^{-1}(dx+idy)$ satisfies 
$d\theta=-2\pi i\omega$. And
\begin{equation}\label{auxr2}
\int_0^1\theta(X_t)dt=-2k\pi i \frac{x_0^2+y_0^2}{x_0^2+y_0^2+1}.
\end{equation}  
From (\ref{actint}), (\ref{auxr1}) and (\ref{auxr2}) it follows 
${\frak A}_{\psi}(q)=k/2 +{\Bbb Z}$
 and 
\begin{equation}\label{kappapsi}
\kappa(\psi)=e^{ik\pi}.
\end{equation}
  
\smallskip

Next we determine the value of $\kappa(\psi)$ by using the results of Section 3.
First of all the prequantum bundle for $({\Bbb CP}^1,\omega)$ is the hyperplane bundle
\cite{pG78}
on ${\Bbb CP}^1$. On the other hand the  functional 
$$ci \in{\frak u}(1)\subset{\frak su}(2)\mapsto 2\pi i\eta(\text{diag}(ci,\,-ci))=-ic$$
is the derivative of $\Lambda:g\in U(1)\mapsto g^{-1}\in U(1)$. Therefore 
the respective prequantum data are $SU(2)$-invariant.
  The isotopy $\{\psi_t\}$ of ${\Bbb CP}^1$
determines the vector fields $X_{A_t}$, where $A_t=\text{diag}(-i\Dot a_t,\,i\Dot a_t)$.
In this case the solution  to $\Dot h_t h_t=A_t$ is  
$h_t=\text{diag}(e^{-ia_t},\, e^{ia_t})$.
Hence, by Theorem \ref{Tmcoador}, 
$$\kappa(\psi)=\Lambda(h_1)=h_1^{-1}=e^{ik\pi}.$$
This result agrees with (\ref{kappapsi}).

\medskip   

{\sc The invariant $\kappa$ of a Hamiltonian flow in $S^2$.}
 For $G=SU(2)$, if
$$\eta:\pmatrix
ai& w \\
-\bar w& -ai
\endpmatrix\in{\frak su}(2)\mapsto\frac{na}{2\pi}\in{\Bbb R},$$
with $n\in{\Bbb Z}$, then the orbit ${\cal O}_{\eta}=SU(2)/U(1)=S^2$ admits 
and $SU(2)$-invariant quantization and the corresponding character 
$\Lambda$ of $U(1)$ is $\Lambda(z)=z^n$.

  Let $\tilde\psi_t$ be the symplectomorphism of 
 $S^2$ given by 
$$\tilde\psi_t(q)=\text{exp}(t(aA+bB))\cdot q,$$
 where $a,b\in{\Bbb R}$ and $A,B$ are the matrices introduced in (\ref{matrices}) . For 
$t_1=(a^2+b^2)^{-1/2}\pi$,
 $\tilde\psi_{t_1}=\text{id}$; in fact
\begin{equation}\label{exponSU}
\text{exp}(t(aA+bB))=\begin{pmatrix} \cos |c| & \epsilon\sin|c| \\
-\bar\epsilon & \cos|c| \end{pmatrix},
\end{equation}
with c=t(b+ai) and $\epsilon=c/|c|$
 (see \cite{aV01}). If we set 
$$E:=\pi(a^2+b^2)^{-1/2}(aA+bB),$$
 by (\ref{exponSU})
$\text{exp}(E)=-\text{Id}$. So  the family
$\{\psi_t\}_{t\in[0,1]}$, defined by $\psi_t(q)=\text{exp}(tE)q$, is a closed 
Hamiltonian flow
on the orbit ${\cal O}_{\eta}$. By Theorem \ref{Tmcoador}
$$\kappa(\psi)=\Lambda(e^E)=\Lambda(-\text{Id})=(-1)^n.$$
This result agrees with that one obtained in \cite[Theorem 21]{aV01} by direct calculation.

This result can be deduced from (\ref{Weylf}), when $n<0$. Here Lax equation $\Dot h_t h_t^{-1}=E$
has the solution $h_t=\text{exp}(tE)$. 
The Weyl's character formula \cite{JJ00} is very simple for
 the group $SU(2)$; in this case, there is
only one positive root $\alpha$ and
the Weyl group has only two elements. 
We take for $\alpha$ the linear map defined by
$$\alpha(\text{diag}(ai,\,-ai))=2ai;$$ 
so $\alpha^{\vee}=\text{diag}(1,\,-1)$.
As we assume that $n<0$, then $-\lambda:=-2\pi i\eta$ is the 
highest weight of a representation $\pi$ of $SU(2)$.
 For $t\in U(1)$, $t^{\lambda}=t^{-n}$ and $t^{\alpha}=t^2$. Therefore (see \cite{JJ00})
$$\text{dim}\, \pi=-n+1 \,\,\,\text{and}\,\,\, \chi_{\pi}(t)=\sum_{k=0}^{-n}t^{-n-2k}.$$
 Hence 
$$\chi_{\pi^*}(h_1)=\chi_{\pi}(-1)=(-n+1)(-1)^n,$$
 and from Corollary \ref{Corchi} we again obtain  the value $(-1)^n$ for $\kappa(\psi)$.


\bigskip 




\begin{thebibliography}{99}   
 
  
   
 \bibitem{JJ00}
 J.J. Duistermaat, J.A.C. Kolk  
 {\em Lie Groups.} 
  Springer, 
   Berlin 
  (2000) 
 

  
\bibitem{wF91}
W. Fulton,   J.  Harris,  
 {\em Representation Theory.}
Springer-Verlag,
New York (1991)



 
 \bibitem{pG78}
  P. Griffiths,   J.  Harris,  
 {\em Principles of algebraic geometry.} 
  John Wiley and Sons, 
  New York 
  (1978) 
 
 
    \bibitem{vG82}
  V. Guillemin,   S.  Sternberg,  
 {\em Geometric quantization and multiplicites.}
   Invent. Math. 
{\bf 67}
   (1982)
  515-538 
 
  


   
  


\bibitem{aK76}
A.A. Kirilov
{\em  Elements of the Theory of Representations.} 
  Springer-Verlag, 
  Berlin, Heidelberg 
 (1976)

  
 



\bibitem{KN}
 S. Kobayashi,  K.  Nomizu,  
 {\em Foundations of differential geometry I.}
  Wiley, 
  New York 
 (1963) 


  
     
   
\bibitem{bK70}
 B. Kostant
{\em Quantization and unitary representations.} In Lectures in modern analysis III
(ed. C.T. Taam). Lecture notes in Mathematics, vol. 170. 87-208. Springer-Verlag, Berlin (1970)

    
 
 
\bibitem{Mc-S}
  D. McDuff, D. Salamon {\em Introduction to symplectic topology.}
Clarenton Press,
 Oxford 
(1998)

 

\bibitem{jS80}
J. Sniatycki,
{\em Geometric Quantization and Quantum Mechanics.}
Springer-Verlag,
New York  (1980)



  
 
  


\bibitem{aV01}
    A. Vi\~{n}a,
 {\em  Hamiltonian symplectomorphisms and the Berry phase.} 
  (To appear in J. Geom. Phys.)  (math.SG/0009206)





 \bibitem{aW89}    A. Weinstein,
 {\em Cohomology of Symplectomorphism Groups and Critical Values of Hamiltonians.} 
    Math. Z.  
 {\bf 201} 
 (1989) 
 75--82




 
 \bibitem{aW90}
    A. Weinstein,
 {\em Connections of Berry and Hannay type for moving Lagrangian Submanifolds.} 
    Advances in Mathematics 
 {\bf 82} 
 (1990) 
 133--159

   
  \bibitem{nW92}
 N.M.J. Woodhouse, 
 {\em Geometric quantization.}
  Clarenton Press,    Oxford 
  (1992)

 
\end{thebibliography}
\end{document}